\documentclass[a4paper, BCOR=0.0mm, DIV=calc]{scrartcl}
\usepackage[T1]{fontenc}
\usepackage[utf8]{inputenc}
\usepackage[ngerman]{babel}
\usepackage[osf,sc]{mathpazo}
\usepackage{textcomp}
\usepackage{microtype}
\usepackage{amsmath, amssymb, amsfonts}
\usepackage{tabularx, booktabs}
\KOMAoptions{twoside=false, twocolumn=false, headinclude=false, footinclude=false, mpinclude=false, pagesize=auto}
\recalctypearea

\newcommand{\etc}{\text{etc}}
 
\setkomafont{section}{\mdseries\scshape\Large}
\setkomafont{title}{\mdseries\scshape}

\begin{document}
\title{Zahlentheoretische Theoreme, mit einer neuen Mehthode bewiesen\footnote{
Originaltitel: "`Theoremata arithmetica nova methodo demonstrata"', erstmals publiziert in "`\textit{Novi Commentarii academiae scientiarum Petropolitanae} 8, 1763, pp. 74-104"', Nachdruck in "`\textit{Opera Omnia}: Series 1, Volume 2, pp. 531 - 555"' und "`\textit{Commentat. arithm.} 1, 1849, pp. 274-286 [E271a]"', Eneström-Nummer E271, übersetzt von: Alexander Aycock, Textsatz: Artur Diener,  im Rahmen des Projektes "`Euler-Kreis Mainz"' }}
\author{Leonhard Euler}
\date{}
\maketitle
\subsection*{Einleitung}
Außer den verschiedenen Rechenarten, die im allgemeinen in der Zahlentheorie zu behandelt werden pflegen und quasi den praktischen Teil dieser Disziplin festsetzen, ist der theoretische Teil derselben, der sich mit dem Untersuchen der Natur der Zahlen beschäftigt, schon einst begonnen worden nicht weniger behandelt zu werden, wie sich aus Euklid und Diophant einsehen lässt, wo man außerordentliche Eigenschaften der Zahlen entdeckt und bewiesen vorfindet. Je mehr aber darauf die Mathematiker die Gestalt und die Beschaffenheit der Zahlen untersucht haben, haben sie um Vieles mehr Eigenschaften derer beobachtet, woher sie die schönsten Theoreme, die die Natur der Zahlen aufzeigen, berechnet haben, die teils durch Beweise untermauert worden sind, teils diese immer noch brauchen, weil sie entweder von den Autoren nicht gefunden worden sind oder mit der Zeit verloren gegangen sind; auf diese Art tauchen verstreut viele zahlentheoretische Theoreme dieser Art auf, deren Beweise man noch heute ersehnt, auch wenn sich deren Gültigkeit nicht in Zweifel ziehen lässt. Und hier müssen wir einen riesigen Unterschied, welcher zwischen zahlentheoretischen und geometrischen Theoremen einhergeht, nicht wenig bewundern, weil kaum irgendeine Proposition hervorgebracht werden kann, die nicht klar ist, entweder Wahres oder Falsches zu zeigen, während dagegen viele Propositionen über die Natur der Zahlentheorie bekannt sind, deren Gültigkeit sich von uns erkennen lässt, aber sich in kleinster Weise beweisen lässt. Man hat eine große Menge von Theoremen dieser Art, die von Fermat hinterlassen wurde, deren Beweise er versicherte zum größten Teil gefunden zu haben, welche zum außerordentlichen Verlust für diese Wissenschaft zu bedauern sind, mit seinen verloren gegangen zu sein. Wie viele Beweise solcher Theoreme aber entweder bekannt sind oder wiederentdeckt sind, wird sich bei diesen gewiss um vieles größere Kraft des Geistes zeigen, die wir kaum bei irgendeiner anderen Art von Beweisen entdecken; daher ist bei dieser Aufgabe nicht so sehr die Nützlichkeit, durch die Wissenschaft der Zahlen illustriert wird, zu schätzen, wie die größte Freiheit, durch die die Beweise dieser Art sich in Bezug auf andere unterscheiden. Und deswegen, weil ich schon öfter, als es den Meisten als angemessen erscheinen kann, in diesem Gebiet gearbeitet habe, glaube ich für meine Person, die Arbeit nicht verschwendet zu haben und glaube auch jetzt nicht, dass die Theoreme, die ich hier vorlege, frei von Nutzen sind. Besonders bemerkenswert schien jenes Theorem von Fermat, in welchem er versichert, dass in dieser Formel $a^{p-1} - 1$ enthaltenen Zahlen immer durch die Zahl $p$ teilbar sind, wenn sie natürlich prim war und auch $a$ trotzdem durch sie keine Teilung zulässt, von welchem Theorem ich schon zwei Beweise gegeben habe. Nun betrachte ich selbiges, aber in weiterem Sinne und untersuche es in der Art, wenn der Teiler keine Primzahl war, sondern irgendeine Zahl $N$, ein Exponent von welcher Art der Potenz zugeteilt werden muss, dass der Ausdruck $a^n - 1$ immer durch die Zahl $N$  teilbar ist, während die Zahl $a$ mit ihr keinen gemeinsamen Teiler hat. Ich habe aber entdeckt, dass dies immer passiert, sooft der Exponent $n$ gleich der Menge der Teiler der Zahlen kleiner als $N$ war, die zu $N$ prim sind. Um das also zu beweisen, braucht man vor allem Theoreme solcher Art, aus denen für irgendeine vorgelegte Zahl $N$ erkannt werden kann, wie viele unter den Zahlen kleiner als sie selbst zu ihr prim sein werden, oder die mit ihr keinen gemeinsamen Teiler haben; diese Theoreme scheinen nun selbst einen weiteren Nutzen zu haben und zu anderen verborgeneren Eigenschaften der Zahlen einen Zugang zu verschaffen. Nachdem diese Dinge vorausgeschickt worden sind, ist der Beweis der vorgelegten Wahrheit deart beschaffen, dass er größerer Aufmerksamkeit nicht unwürdig sein wird.
\section*{Theorem 1}
\paragraph{§1}
Wenn die Terme einer arithmetischen Progression, deren Differenz eine zu $n$ prime Zahl sei, durch irgendeine Zahl $n$ geteilt wurden, werden unter den Resten alle Zahlen kleiner als die Zahl $n$ auftauchen.
\subsection*{Beweis}
Es sei der erste Term der arithmetischen Progression gleich $a$ und die Differenz gleich $d$, die eine zu $n$ prime Zahl sei oder die mit der Zahl $n$ keinen gemeinsamen Teiler außer der Einheit habe, sodass die arithmetische Progression
\[
	a,\quad a+d,\quad a+2d,\quad a+3d,\quad a+4d,\quad a+5d,\quad \etc
\] 
sein wird, und ich sage, wenn die einzelnen Terme durch die Zahl $n$ geteilt werden, dass unter den Resten alle Zahlen kleiner als $n$ auftauchten. Um das zu zeigen, wird es genügen, nur die $n$ Terme dieser Progression betrachtet zu haben, die
\[
	a,\quad a+d,\quad a+2d,\quad a+3d,\quad \dots \quad a+(n-1)d
\]
sind. Wenn daher also diese einzelnen Terme durch $n$ geteilt werden, müssen alle Reste zueinander verschieden sein. Wenn nämlich zwei Terme, z.\,B. $a + \mu d$ und $a + \nu d$, während $\mu$ und $\nu$ Zahlen kleiner als $n$ selbst sind, diese durch $n$ geteilt gleiche Reste liefern würden, würde deren Differenz $(\nu - \mu )d$ jedenfalls durch $n$ teilbar sein. Weil aber die Zahlen $d$ und $n$ keinen gemeinsamen Teiler haben, ist es notwendig, dass $(\nu - \mu )$ eine Teilung durch $n$ zuließen. Weil daher all jene Reste verschieden sind und ja die Anzahl der Terme von der Zahl her gleich $n$ ist, werden bei diesen natürlich alle Zahlen kleiner als $n$ auftauchten, natürlich 
\[
	0,\quad 1,\quad 2,\quad 3,\quad 4,\quad 5,\quad \dots\quad  n-1,
\]
wenn freilich die Differenz der Progression $d$ eine zum vorgelegten Teiler $n$ prime Zahl ist.

\hfill\textsc{Q.E.D.}
\subsection*{Korollar 1}
\paragraph{§2}
Unter diesen Termen dieser arithmetischen Progression, deren Anzahl $n$ ist, solange seine Differenz eine zu $n$ prime Zahl ist, wird man also gewiss eine finden, die durch $n$ teilbar ist; dann aber wird auch eine da sein, die durch $n$ geteilt einen gegebenen Rest $r$ zurücklassen wird.
\subsection*{Korollar 2}
\paragraph*{§3}
Wenn also die Zahl $d$ zu $n$ prim war, kann immer eine Zahl dieser Zahl dieser Form $a + \nu d$ beschafft werden, während $a$ irgendeine Zahl ist und $\nu$ kleiner als $n$ ist, die durch die Zahl $n$ teilbar ist, und es wird auch unter denselben Bedingungen immer eine solche Zahl $a+rd$ gegeben sein, die durch $n$ geteilt einen gegebenen Rest $r$ lassen wird.
\subsection*{Korollar 3}
\paragraph{§4}
Nachdem also Zahlen $a$ und $d$ gegeben worden sind, von denen diese $d$ zu $n$ prim sei, lassen sich immer Zahlen $\mu$ und $\nu$ finden, dass dieser Gleichung
\[
	a + \nu d = \mu n
\]
oder auch dieser
\[
	a + rd = \mu n + r
\]
genügt wird, was für eine Zahl kleiner als $n$ für $r$ auch immer angenommen wird.
\subsection*{Bemerkung}
\paragraph{§5}
Was wir über die Anzahl $n$ der Terme der arithmetischen Progression bewiesen haben, gilt auch über die ganze Progression ins Unendliche fortgesetzt; die Terme, die nach jenen $n$ Termen folgen, erzeugen nämlich in derselben Reihenfolge die Reste, wenn sie durch $n$ geteilt werden. So stimmen die Reste der Terme, die nach $a + (n-1)d$ folgen und $a+nd$, $a+(n+1)d$, $a+(n+2)d$, etc sind, durch $n$ geteilt mit den Resten überein, die aus den anfänglichen Termen $a$, $a+d$, $a+2d$, etc entstehen. Und wenn die ganze Reihe in unendlich viele Perioden aufgeteilt werden, indem man irgendeiner $n$ Terme auf diese Weise zuteilt
\[
	a, \quad a+d, \dots, a+(n-1)d | a+nd, \dots, a+(2n-1)d | a+2nd, \dots, a+(3n-1)d | \dots
\] 
werden die Terme einer beliebigen Periode dieselben Reste, in derselben Reihenfolge geordnet, liefern; die ersten und zweiten und dritten Terme aller Perioden werden nämlich gleich bleiben die gleichen Reste geben. Wenn wir daher die Art der Reste erkennen wollen, genügt es, eine einzige Periode untersucht zu haben.
\section*{Theorem 2}
\paragraph{§6}
In einer arithmetischen Progression, deren Anzahl an Termen gleich $n$ ist, werden so viele Terme zur Zahl $n$ prim sein, wie unter den Zahlen kleiner als $n$ selbst prime zu $n$ gegeben sind, solange die Differenz der Progression zu $n$ prim war.
\subsection*{Beweis}
Es sei nämlich $a$ der erste Term und $d$ die Differenz der Progression, die zu $n$ prim sei, und daher ist die Progression, die $n$ Terme enthält,
\[
	a,\quad a+d,\quad a+2d,\quad a+3d,\quad \dots \quad a +(n-1)d. 
\]
Weil ja also, wenn diese Terme durch $n$ geteilt werden, natürlich alle Reste der Zahl kleiner als $n$ selbst auftauchen, wollen wir setzen, dass aus irgendeinem Term $a+\nu d$ der Rest $r$ resultiert und es ist klar, wenn $r$ eine zu $n$ prime Zahl war, dass auch jener Term $a + \nu d$ zu $n$ prim sein wird; wenn aber $r$ mit $n$ einen gemeinsamen Teiler hat, wird selbiger auch gemeinsamer Teiler der Zahlen $n$ und $a+\nu d$ sein. Wie viele Zahlen daher unter den Zahlen kleiner als $n$ prim zu $n$ waren, man wird genauso viele zu $n$ prime Zahlen unter den Termen der vorgelegten arithmetischen Progression haben.

\hfill\textsc{Q.E.D.}
\subsection*{Korollar 1}
\paragraph{§7}
Wenn $n$ eine Primzahl war, weil alle Zahlen kleiner al selbige zu dieser auch prim sind, deren Anzahl also gleich $n-1$ ist, werden auch in jener arithmetischen Progression alle Terme aus einem einzigen zu $n$ prim sein, weil ja ein einziger durch $n$ teilbar ist.
\subsection*{Korollar 2}
\paragraph{§8}
Wenn aber $n$ eine zusammengesetzte Zahl war, sind zwischen den Zahlen kleiner als selbige solche gegeben, die mit ihr einen gemeinsamen Teiler haben, und man wird in der Tat auch so viele in der arithmetischen Progression finden, bei denen dieselben gemeinsamen Teiler mit $n$ übereinstimmen.
\subsection*{Korollar 3}
\paragraph{§9}
Wenn $n=6$ war, werden, weil unter den Zahlen kleiner als $6$ zwei zu $6$ prim aind, natürlich $1$ und $5$, in der ganzen arithmetischen Progression der $6$ Terme
\[
	a,\quad a+d,\quad a+2d,\quad a+3d,\quad a+4d,\quad a+5d
\]
nur zwei zu $6$ prim sein, solange die Differenz $d$ eine zu $6$ prime Zahl ist. Wenn man so $a=4$, $d=5$ nimmt, sind zwei dieser sechs Zahlen $4$, $9$, $14$, $19$, $24$, $29$, natürlich $19$ und $29$, zu $6$ prim, eine, $24$, ist durch $6$ teilbar, die übrigen $4$, $9$, $14$ sind in der Tat genauso wie $2$, $3$, $4$ zu $6$ nicht prim.
\subsection*{Bemerkung}
\paragraph{§10}
Diese Theoreme haben in der Lehre und der Betrachtung der Natur der Zahlen einen riesigen Nutzen, hier scheint es aber nur ratsam diese zu verwenden, um diese Frage zu erörtern, nachdem irgendeine Zahl $n$ vorgelegt worden ist, wie viele unter den Zahlen kleiner als $n$ zu derselben Zahl $n$ prim sind. Es ist freilich sofort klar, wenn $n$ eine Primzahl ist, dass alle Zahlen kleiner als selbige zugleich zu ihr prim sein werden und deren Anzahl deshalb gleich $n-1$ sein wird. Wenn aber $n$ eine zusammengesetzte Zahl ist, ist die Menge der zu ihr primen Zahlen kleiner als selbige; wie groß sie aber in jedem Fall ist, kann nicht so leicht angegeben. Wenn $n=12$ ist, findet man unter den Zahlen kleiner als diese nur vier zu $12$ prime, natürlich $1$, $5$, $7$, $11$, und wenn $n=60$ ist, sind die kleineren zu dieser Zahl primen 
\[
	1,~~ 7,~~ 11,~~ 13,~~ 17,~~ 19,~~ 23,~~ 29,~~ 31,~~ 37,~~ 41,~~ 43,~~ 47,~~ 49,~~ 53,~~ 59 
\] 
deren Anzahl $16$ ist, woher die übrigen $43$ alle mit $60$ einen gemeinsamen Teiler haben. Hier sollte man sich aber daran erinnern, dass die Einheit eine zu allen Zahlen prime Zahl ist, auch wenn sie ein Teiler aller ist; das ist aus der Definition klar, nach welcher Zahlen bezeichnet werden, zueinander prim zu sein, die außer der Einheit keinen anderen Teiler haben.  
\section*{Theorem 3}
\paragraph{§11}
Wenn $n$ irgendeine Potenz der Primzahl $p$ oder $n=p^m$ ist, werden zwischen den Zahlen größer als selbige so viele zu ihr prim sein, wie Einheiten in
\[
	p^m - p^{m-1} = p^{m-1}(p-1)
\]
enthalten sind.
\subsection*{Beweis}
Die Menge aller Zahlen kleiner als die Potenz $n=p^m$ ist $p^m - 1$, unter diesen findet man aber einige, die zu $n$ nicht prim sind, natürlich alle Vielfachen von $p$ kleiner als $n$ und zusätzlich keine anderen; daraus werden die folgenden Zahlen zu $n$ nicht prim sein
\[
	p,\quad 2p,\quad 3p,\quad 4p, \quad \dots \quad p^m - p
\]
deren Anzahl $p^{m-1}-1$ ist; nachdem diese von der Zahl aller $n=p^m$ abgezogen wurde, findet man die Menge der kleineren, die zu $p^m$ prim sind, deren Anzahl deshalb gleich $p^m - p^{m-1} = p^{m-1}(p-1)$ ist.

\hfill\textsc{Q.E.D}.
\subsection*{Korollar 1}
\paragraph{§12}
Daher folgt also, was an sich klar ist, wenn $n=p$ war, während $p$ eine Primzahl ist, dass die Anzahl aller zu ihr primen Zahlen kleiner als selbige gleich $p-1$ ist, weil ja alle Zahlen kleiner als dieselbe zugleich zu ihr prim sind.
\subsection*{Korollar 2}
\paragraph{§13}
Aber wenn $n=p^2$ ist, es die Menge unter den Zahlen kleiner als selbige von diesen, die zu ihr prim sind, gleich $pp-p = p(p-1)$; die übrigen, deren Anzahl $p-1$ ist, werden zu $n=p^2$ nicht prim sein oder durch $p$ teilbar.
\subsection*{Korollar 3}
\paragraph{§14}
Nachdem aber irgendeine Potenz $n=p^m$ der Primzahl vorgelegt wurde, findet man unter den Zahlen kleiner als selbige, deren Menge gleich $p^m - 1$ ist, $p^{m-1}-1$, die durch $p$ teilbar sind und daher zu $p^m$ nicht prim sind; aber alle übrigen, deren Anzahl gleich $p^m - p^{m-1} = p^{m-1}(p-1)$ ist, sind zu $p^m$ prim.
\subsection*{Bemerkung}
\paragraph{§15}
Wenn also die vorgelegte Zahl $n$ eine Potenz einer Primzahl $n$ war, werden wir mithilfe dieser Regel angeben können, wie viele unter allen Zahlen kleiner als selbige zu dieser prim sein werden. Wann immer aber die Zahl $n$ aus zweien oder mehr Primzahlen zusammengesetzt war, kann daher diese Frage nicht erledigt werden; durch Anwenden der vorhergehenden Theoreme können wir diese sich weiter erstreckende Frage lösen.
\section*{Theorem 4}
\paragraph{§16}
Wenn die Zahl $n$ das Produkt zweier Primzahlen $p$ und $q$ war oder $n=pq$, ist die Menge aller zu ihr primen Zahlen kleiner als selbige gleich $(p-1)(q-1)$.
\subsection*{Beweis}
Weil die Anzahl aller Zahlen kleiner als $n=pq$ gleich $pq-1$ ist, müssen daher die zuerst ausgeschlossen werden, die durch $p$ teilbar sind, darauf aber auch die, die durch $q$ teilbar sind; und nach dem Streichen dieser wird man die gesuchte Menge finden. Man notiere also von der Einheit bis hin zu $pq$ die Zahlen, die zu $p$ prim sind, auf diese Weise
\[
\begin{array}{cccccc}
1 & 2 & 3 & 4 & \dots & p-1 \\
p+1 & p+2 & p+3 & p+4 & \dots & 2p-1 \\
2p+1 & 2p+2 & 2p+3 & 2p+4 & \dots & 3p-1 \\
3p+1 & 3p+2 & 3p+3 & 3p+4 & \dots & 4p-1 \\
\vdots & \vdots & \vdots & \vdots &  & \vdots \\
(q-1)p+1 & (q-1)p+2 & (q-1)p+3 & (q-1)p+4 & \dots & pq-1
\end{array}
\]
und nun müssen aus diesen nur die ausgewählt werden, die zugleich auch zu $q$ prim sind. Man betrachte also die vertikalen Reihen, deren Anzahl $p-1$ ist; eine beliebige aber enthält $q$ Terme, die in einer arithmetischen Progression wachsen, während die Differenz $p$ ist, die eine zu $q$ prime Zahl ist. In einer beliebigen vertikalen Reihe werden also alle Terme außer einem zu $q$ prim sein (nach §$7$); daher enthält jede vertikale Reihe $q-1$ zu $q$  prime Zahlen. Weil daher die Anzahl der vertikalen Reihen $p-1$ ist, sind in allen gleichzeitig $(p-1)(q-1)$ zur Zahl $q$ prime Zahlen enthalten und dieselben werden also auch zum Produkt $pq$ prim sein; als logische Konsequenz wird man unter allen Zahlen kleiner als $pq$ selbst $(p-1)(q-1)$ zu $pq$ prime Zahlen finden.

\hfill\textsc{Q.E.D.}
\subsection*{Korollar 1}
\paragraph{§17}
Weil die Menge aller Zahlen kleiner als das Produkt $pq$ gleich $pq-1$ ist, sind unterdessen immer $(p-1)(q-1) = pq-p-q+1$ zu $pq$ prim, die übrigen aber, deren Anzahl $p+q-2$ ist, sind zu ihr nicht prim oder haben mit ihr als gemeinsamen Teiler entweder $p$ oder $q$.
\subsection*{Korollar 2}
\paragraph{§18}
Das ist auch daher klar, weil unter den Zahlen kleiner als das Produkt $pq$ $q-1$ durch $p$ teilbare Zahlen sind, natürlich
\[
	p,\quad 2p,\quad 3p,\quad 4p,\quad \dots \quad (q-1)p
\]
darauf unter denselben $p-1$ Zahlen durch $q$ teilbar sind, nämlich
\[
	q,\quad 2q,\quad 3q,\quad 4q,\quad \dots \quad (p-1)q;
\]
weil ja von jenen alle verschieden sind, wird man im Ganzen
\[
	(q-1) + (p-1) = p+q-2
\]
Zahlen haben, die zu $pq$ nicht prim sind.
\subsection*{Korollar 3}
\paragraph{§19}
Wenn man also sucht, wie viele von den Zahlen $1$ bis hin zu $15$ zu $15$ prime Zahlen sind, lehrt wegen $p=3$ und $q=5$ die Regel, dass die Anzahl derer $2\cdots 4 = 8$ ist, die ja sind
\[
	1, \quad 2,\quad 4,\quad 7,\quad 8,\quad 11,\quad 13,\quad 14
\]
sind. Auf die gleiche Weise ist die Menge der zu $35$ primen Zahlen von $1$ bis $35$ wegen $p=5$ und $q=7$ $4\cdot 6 = 24$ und diese Zahlen sind
\[
	1,\,2,\,3,\,4,\,5,\,6,\,8,\,9,\,11,\,12,\,13,\,16,\,17,\,18,\,19,\,22,\,23,\,24,\,26,\,27,\,29,\,31,\,32,\,33,\,34
\]
\subsection*{Bemerkung}
\paragraph{§20}
Weil ja hier die Frage über Zahlen geht, die zu einer Zahl prim sind und kleiner als sie, werden diese sich angenehm als zu dieser Zahl prime Teile bezeichnen lassen. Wenn so die Primzahl $p$ vorgelegt war, wird die Anzahl der zu ihr primen Teile gleich $p-1$ sein; wenn die vorgelegte Zahl eine Potenz der Primzahl, $p^n$, ist, wird die Anzahl der zu ihr primen Zahlen gleich $p^n - p^{n-1}(p-1)$ sein; aber wenn die vorgelegte Zahl das Produkt zweier verschiedenen Primzahlen gleich $pq$ ist, ist die Anzahl der zur ihr primen Teile gleich $(p-1)(q-1)$; und auf diese Weise können wir die Umschweife beim Sprechen verkürzen; auf die gleiche Weise können wir zeigen, wenn die vorgelegte Zahl das Produkt aus drei verschiedenen Primzahlen gleich $pqr$ ist, dass die Anzahl der zu ihr primen Teile gleich $(p-1)(q-1)(r-1)$ sein wird; und das ließe sich auch auf Produkte mehrerer ausdehnen. Aber die folgende Proposition umfasst all diese Fälle in sich.
\section*{Theorem 5}
\paragraph*{§21}
Wenn $A$ und $B$ zwei zueinander prime Zahlen sind und die Anzahl der zu $A$ primen Teile gleich $a$ ist, die Anzahl der zu $B$ primen Teile gleich $b$ ist, dann wird die Anzahl der zum Produkt $AB$ primen Teile gleich $ab$ sein.
\subsection*{Beweis}
Es seien $1$, $\alpha$, $\beta$, $\gamma$, $\dots$, $\omega$ zu $A$ prime Zahlen und kleiner als selbige oder zu $A$ prime Teile, die Anzahl welcher Teile also per Annahme gleich $a$ ist. Es werden also so viele Zahlen zu $A$, wie von $A$ zu $2A$ prim sein, und ebenso von $2A$ zu $3A$, usw. Auf diese Weise können alle zu $A$ prime Zahlen von der Einheit bis hin zur vorgelegten Zahl $AB$ beschafft werden, welche das folgende Schema beschaffen wird:
\[
\begin{array}{ccccc}
1 & \alpha & \beta & \dots & \omega \\
A+1 & A+\alpha & A+\beta & \dots & A+\omega \\
2A + 1 & 2A+\alpha  &2A+\beta & \dots & 2A + \omega \\
3A+1 & 3A+\alpha & 3A+\beta & \dots & 3A + \omega \\
\vdots & \vdots & \vdots & & \vdots \\
(B-1)A+1 & (B-1)A+\alpha & (B-1)A+\beta & \dots &  (B-1)A+\omega 
\end{array}
\]
Hier enthalten die einzelnen horizontalen Reihen $a$ Terme und die Anzahl aller horizontalen Reihen ist gleich $B$, woher alle Reihen zusammen $aB$ Terme ergeben, die schon alle zu $A$ prim sein werden. Daher müssen also noch die ausgeschlossen werden, die zu $B$ nicht prim sind, damit auf diese Weise die zurückgelassen werden, die nicht nur zu $A$, sondern auch zu $B$ und daher zum Produkt $AB$ prim sind; oder es müssen aus diesen Reihen nur die Terme gezählt werden, die auch zu $B$ prim sind. Für dieses Ziel wollen wir die Reihen vertikal betrachten; und weil die Anzahl der vertikalen Reihen gleich $a$ ist, wird eine beliebige vertikale Reihe $B$ Terme enthalten, die in einer aritmetischen Progression vermehrt worden sind; weil deren Differenz gleich $A$ ist und eine zu $B$ prime Zahl ist, wird durch Theorem $2$ eine beliebige vertikale Reihe so viele zu $B$ prime Terme enthalten, wie zur Zahl $B$ prime Teile gegeben sind. Die Anzahl derer ist also per Annahme gleich $b$. Weil also die einzelnen vertikalen Reihen $b$ zu $B$ prime Terme enthalten, die deshalb auch zum Produkt $B$ prim sein werden, wird die Anzahl aller zu $AB$ primen Terme, das heißt zu dieser Zahl $AB$ primen Teile, gleich $ab$ sein.

\hfill\textsc{Q.E.D.} 
\subsection*{Korollar 1}
\paragraph*{§22}
Wenn darüber heraus eine dritte Zahl $C$ hinzugefügt wird, die zu jeder der beiden vorhergehenden $A$ und $B$ oder zum Produkt der $AB$ prim ist, und die Anzahl der zu $C$ primen Teile gleich $c$ ist, dann wird die Anzahl der zum Produkt $ABC$ primen Teile gleich $abc$ sein. Es kann nämlich das Produkt $AB$ als eine Zahl betrachtet werden, die Menge der zu ihr primen Teile gleich $ab$ ist; und weil $c$ zu $AB$ prim ist, hat das Theorem hier Geltung.
\subsection*{Korollar 2}
\paragraph*{§23}
Weil also jede einzelne Zahl in $N$ zueinander prime Faktoren aufgelöst werden kann, welche einzelnen entweder selbst Primzahlen sind oder Potenzen von Primzahlen, wird mit Hilfe dieser Regel die Menge der zu irgendeiner Zahl $N$ primen Fälle angegeben werden können.
\subsection*{Korollar 3}
\paragraph*{§24}
Während natürlich $p$, $a$, $r$, $s$ etc Primzahlen sind, wird jeder Zahl $N$ in einer Form dieser Art $N = p^{\lambda} q^{\mu}r^{\nu}s^{\xi}$ erfasst werden, woher die Anzahl der zu $N$ primen Teile
\[
	p^{\lambda-1}(p-1)\cdot p^{\mu - 1}(q-1)\cdot r^{\nu - 1}(\nu - 1)\cdot s^{\xi - 1}(s-1)
\]
sein wird.
\subsection*{Korollar 4}
\paragraph*{§25}
Für einfachere Formen von Zahlen wird sich die Menge der zu diesen primen Teile so verhalten:

\begin{tabular}{p{1.7cm} p{4cm} >{\centering\arraybackslash}p{1.7cm} >{\centering\arraybackslash}p{3cm}}
\toprule
\small\centering Vorgelegte Zahl & \small \centering Menge der zu ihr primen Teile & \small\centering Vorgelegte Zahl & \small \centering Menge der zu ihr primen Teile 
\tabularnewline
\midrule
\small
$p$ & $p-1$ & $2$ & $1$ \\
\midrule
$pp$ & $p(p-1)$ & $3$ & $2$ \\
$pq$ & $(p-1)(q-1)$ & $4$ & $2$ \\
\midrule
$p^3$ & $pp(q-1)$ & $5$ & $4$ \\
$p^2q$ & $p(p-1)(q-1)$ & $6$ & $2$ \\
$pqr$ & $(p-1)(q-1)(r-1)$ & $7$ & $6$ \\
\midrule
$p^4$ & $p^3(p-1)$ & $8$ & $4$ \\
$p^3 q$ & $p^2(p-1)(q-1)$ & $9$ & $6$ \\ 
$p^2q^2$ & $p(p-1)q(q-1)$ & $10$ & $4$ \\
$p^2qr$ & $p(p-1)(q-1)(r-1)$ & $11$ & $10$ \\
$pqrs$ & $(p-1)(q-1)(r-1)(s-1)$ & $12$ & $4$ \\
\midrule
$p^5$ & $p^4(p-1)$ & $13$ & $12$ \\
$p^4q$ & $p^3(p-1)(q-1)$ & $14$ & $6$ \\
$p^3q^2$ & $p^2(p-1)q(q-1)$ & $15$ & $8$ \\
$p^2qr$ & $p^2(p-1)(q-1)(r-1)$ & $16$ & $8$ \\
$p^2q^2r$ & $p(p-1)q(q-1)(r-1)$ & $17$ & $16$ \\
$p^2qrs$ & $p(p-1)(q-1)(r-1)(s-1)$ & $18$ & $6$ \\
$pqrst$ & $(p-1)(q-1)(r-1)(s-1)(t-1)$ & $19$ & $18$ \\
\bottomrule
\end{tabular}
\subsection*{Korollar 5}
\paragraph*{§26}
Nachdem daher also irgendeine Zahl vorgelegt wurde, wird die Menge der zur ihr primen Zahlen angenehm bestimmt werden. Wenn z.\,B. 360 vorgelegt wird, wird, weil $360 = 2^3 \cdot 3^2 \cdot 5$ ist, die Menge der zu 360 primen Teile gleich $4\cdot 9\cdot 4 = 96$ sein.
\subsection*{Bemerkung}
\paragraph*{§27}
Dies kann über die Menge der zu einer Zahl primen Teile für unser Unternehmen ausreichen. Dennoch wird es förderlich sein, über die zu einer Zahlen primen Teile selbst dies bemerkt zu haben: Wenn die vorgelegte Zahl $N$ war und unter den zu ihr primen Teilen die Zahl $\alpha$ auftaucht, wird da selbst auch die Zahl $N-\alpha$ auftauchen, weil ja, während $\alpha$ zu $N$ prim ist, auch $N-\alpha$ zu $N$ prim sein wird. Daher wird es also genügen für eine beliebige Zahl nur die Hälfte ihrer kleineren Teile gefunden zu haben, weil die übrigen deren Komplemente zur Zahl $N$ selbst sind. Auf die gleiche Weise wird, wenn $N$ eine gerade Zahl ist, unter den zu $N$ primen Zahlen auch $\frac{1}{2}N-\alpha$ auftauchen, dann aber auch $\frac{1}{2}N+\alpha$. Wenn so $N$ durch irgendeine Zahl $n$ teilbar ist, werden unter den zu ihr primen Teilen auch diese Zahlen auftauchen
\[
	\frac{1}{n}A\pm\alpha, \quad \frac{2}{n}N\pm\alpha,\quad \dots \quad \frac{n-1}{n}N\pm \alpha \quad \text{und}\quad N-\alpha
\]
und daher werden um Vieles leichter die Teile selbst tatsächlich beschafft werden können.
\section*{Theorem 6}
\paragraph*{§28}
Wenn eine Zahl $x$ prim zu $N$ war, dann werden alle Potenzen von $x$ durch $N$ geteilt Reste lassen, die zur Zahl $N$ prim sein werden.
\subsection*{Beweis}
Weil nämlich $x$ eine zu $N$ prime Zahl ist, werden alle Potenzen von ihr auch zu $N$ prim sein, und wenn sie daher durch $N$ geteilt werden, werden auch die Reste zu $N$ prime Zahlen sein.

\hfill\textsc{Q.E.D.}
\subsection*{Korollar 1}
\paragraph*{§29}
Unter den Resten der Potenzen von $x$, die durch $N$ geteilt werden, tauchen also keine anderen Zahlen auf, wenn diese nicht die zu $N$ prime Teile sind; weil die Anzahl deren für die Gestalt der Zahl $N$ bestimmt ist, existieren unzählige Potenzen von $x$, die durch $N$ geteilt gleiche Reste zurücklassen.
\subsection*{Korollar 2}
\paragraph*{§30}
Unter den Resten, die aus der Teilung der Potenzen von $x$ durch die Zahl $N$ entstehen, wird man aber immer die Einheit finden, deshalb weil zwischen den Potenzen von $x$ auch $x^0 = 1$ gezählt werden muss. Ob aber außer der Einheit auch alle übrigen zu $N$ primen Teile unter den Resten auftauchen oder nicht, werden wir bald sehen.
\subsection*{Korollar 3}
\paragraph*{§31}
Wenn man für $x$ die Einheit nimmt, werden alle Reste die Einheit sein, welche Zahl für $N$ auch immer angenommen worden war. Wenn man darauf $x = N-1$ nimmt, welche Zahl zu $N$ auch prim ist, wird man bei den Resten, die aus der Teilung der Potenzen
\[
	(N-1)^0, \quad (N-1)^1, \quad (N-1)^2, \quad (N-1)^3, etc
\]
entstehen, nur zwei verschiedene finden, natürlich 1 und $N-1$, die ununterbrochen abwechselnd auftauchen.
\subsection*{Korollar 4}
\paragraph*{§32}
Je nachdem wie also die Zahl $x$ von der Art zu $N$ beschaffen war, kann es natürlich geschehen, dass unter den Resten aller Potenzen von $x$ nicht alle zum Teiler $N$ primen Teile auftauchen.
\subsection*{Korollar 5}
\paragraph*{§33}
Wenn also alle zur Zahl $N$ primen Teile $1$, $a$, $b$, $c$, $d$, $e$, \dots sind, deren Anzahl gleich $n$ ist, werden unter den erwähnten Resten entweder all diese Teile auftauchen oder nur gewisse, unter denen man aber immer die Einheit finden wird.
\subsection*{Korollar 6}
\paragraph*{§34}
Wenn daher nicht all jene Teile in den Resten, die aus der Teilung der Potenzen von $x$ durch die Zahl $N$ zurückgelassen werden, auftauchen, werden jene Teile in zwei Klassen aufgeteilt werden, deren eine die Teile enthalten wird, die in den Resten auftauchen, die andere aber die Teile, die nicht in den Resten auftauchen.
\section*{Theorem 7}
\paragraph*{§35}
Wenn die Reihe der Potenzen $x^0$, $x^1$, $x^2$, $x^3$, $x^4$, $x^5$, etc durch die Zahl $N$, die zu $x$ prim sei, geteilt wird, werden bis dorthin verschiedene Reste hervorgehen, bis man zu einer Potenz gelangt, die wiederum die Einheit für den Rest liefert.
\subsection*{Beweis}
Weil ja in der Reihe der Potenzen $1, x, x^2, x^3, x^4$, etc ins Unendliche fortgesetzt die Reste nicht alle verschieden sein können, ist nötig, dass schließlich ein bestimmter von den vorhergehenden wieder auftaucht; und ich sage, dass die Einheit dieser Rest ist, weil er als erster von allen auftaucht. Wenn irgendeiner das verneinen sollte, sei $x^{\mu}$ jene Potenz, deren Rest als erstes in den folgenden aus der Potenz $x^{\mu + \nu}$ wieder auftaucht; weil also die Potenzen $x^{\mu}$ und $x^{\mu + \nu}$ die gleichen Reste liefern werden, wird deren Differenz $x^{\mu + \nu} - x^{\mu} = x^{\mu}(x^{\nu}-1)$ durch die Zahl $N$ teilbar sein. Aber der erste Faktor des Produktes $x^{\mu}(x^{\nu}-1)$ ist eine zu $N$ prime Zahl, also ist notwendig, dass der andere $x^{\nu}-1$ durch $N$ teilbar ist. Daher würde aber die Potenz $x^{\nu}$ durch $N$ geteilt den Rest 1 geben und so wird die Einheit unter den folgenden Resten schneller auftauchen als der Rest einer Potenz $x^{\mu}$, welcher ja per Annahme erst in der höheren Potenz $x^{\mu + \nu}$ wiederkehrt. Daraus ist klar, dass kein Rest nochmal auftauchen kann, wenn nicht zuvor die Einheit dazwischen aufgetaucht sein wird.

\hfill\textsc{Q.E.D.}
\subsection*{Korollar 1}
\paragraph*{§36}
Nachdem die Teilung der Terme der Reihe $1$, $x$, $x^2$, $x^3$, $x^4$, etc durch die $N$, die zu $x$ prim ist, von Anfang verschiedene Reste gegeben hatte, z.\,B. $1$, $\alpha$, $\beta$, $\gamma$, etc, wird schließlich wieder der erste Rest $1$ auftauchen; wenn er daher also aus der Potenz $x^{\nu}$ entsteht, wird die Anzahl der vorhergehenden verschiedenen Reste gleich $\nu$ sein.
\subsection*{Korollar 2}
\paragraph*{§37}
Wann immer aber die Potenz $x^{\nu}$ den Rest $1$ gibt, ist es derselbe, den der erste Term $x^0$ gibt, wird die folgende Potenz $x^{\nu + 1}$ denselben Rest geben, den $x^n$ gibt. Und irgendeine der folgenden gibt denselben, welchen die Potenz $x^{\mu}$ gibt. Weil nämlich die Differenz $x^{\nu + \mu}-x^{\mu} = x^{\mu}(x^{\nu}-1)$ durch $N$ teilbar ist, ist notwendig, dass beide Terme $x^{\nu + \mu}$ und $x^{\mu}$ durch $N$ geteilt denselben Rest geben.
\subsection*{Korollar 3}
\paragraph*{§38}
Weil nach der Potenz $x^{\nu}$ dieselben Reste $1$, $\alpha$, $\beta$, $\gamma$, etc der Reihe nach auftauchen, werden die Potenzen $x^{3\nu}$, $x^{4\nu}$, $x^{5\nu}$, etc alle durch $N$ geteilt denselben Rest übrig lassen. Ja es werden sogar alle Potenzen $x^{\mu}$, $x^{\mu + \nu}$, $x^{\mu + 2\nu}$, $x^{\mu + 3\nu}$, $x^{{\mu+4\nu}}$, etc die gleichen Reste liefern.
\subsection*{Korollar 4}
\paragraph*{§39}
Wenn also $x^{\nu}$ die unterste Potenz war, die nach $x^0 - 1$ wiederum die Einheit für den Rest liefert, wird die Anzahl der verschiedenen Reste $\nu$ sein. Weil also die Anzahl der zur Zahl $N$ primen Teile gleich $n$ ist, kann es gewiss nicht geschehen, dass $\nu > n$ ist; es wird also $\nu = n$ oder $\nu < n$ sein.
\subsection*{Korollar 5}
\paragraph*{§40}
Wenn also die Reihe der Potenzen $1$, $x$, $x^2$, $x^3$, etc bis hin zu $x^n$ fortgesetzt wird, wird man unter diesen gewiss zumindest eine einzige außer dem ersten Term $1$ finden, die durch $N$ geteilt die Einheit zurücklässt. Es werden vielleicht irgendwann viele Potenzen dieser Art, aber niemals weniger als eine existieren.
\subsection*{Bemerkung}
\paragraph*{§41}
Die Reste werden ausschließlich immer Zahl kleiner als der Teiler $N$ sein, aber nichts hindert daran, dass wir auch größere Zahlen als Reste betrachten, von welcher Art sie zurückgelassen werden, wenn der Quotient zu klein angenommen wird. Wenn so bei der Teilung der Zahl durch $N$ $N+\alpha$ zurückgelassen wird, muss dieser Rest $\alpha$ selbst äquivalent angesehen werden; und daher, wenn von Resten die Rede ist, sind all diese Zahlen $\alpha$, $N+\alpha$, $2N+\alpha$, $3N+\alpha$ etc gleich einem einzigen Rest $\alpha$ anzusehen. Natürlich verändern irgendwelche Vielfachen des Teilers $N$ entweder hinzugefügt oder von einem Rest $\alpha$ abgezogen seine Natur nicht und auf diese Weise werden auch die negativen Zahlen angenehm unter die Reste gezählt, wie z.\,B. $\alpha - N$ für denselben Rest zu halten ist wie $\alpha$ und der Rest $-1$ dem Rest $N-1$ äquivalent ist. Aus diesen erreicht man, dass alle Zahlen, die durch $N$ geteilt denselben Rest $\alpha$ beschaffen, für denselben Rest gehalten werden können; daraus ensteht nämlich aus einer Zahl durch Teilung ein zu kleiner Quotient, indem man entweder $N+\alpha$ oder $2N+\alpha$ oder $3N+\alpha$ etc, nimmt, aus derselben entsteht, indem man den vollen Quotienten nimmt, der Rest $\alpha$; dann aber wird man eben daher, wenn der Quotient zu groß genommen wird, negative Reste $\alpha - N$ oder $\alpha -2N$ oder $\alpha - 3N$ etc erhalten, die also auch so zu verstehen sind, sich nicht von $\alpha$ zu unterscheiden.
\section*{Theorem 8}
\paragraph*{§42}
Wenn, während die Terme der Progression $1$, $x$, $x^2$, $x^3$, $x^4$, etc durch die zu $x$ prime Zahl $N$ geteilt werden, die Reste $1$, $a$, $b$, $c$ etc waren, werden irgendwelche Produkte entweder zweier oder dreier oder beliebig vieler miteinander multiplizierten auftauchen.
\subsection*{Beweis}
Es mögen also die Reste $a$, $b$, $c$, etc aus den Potenzen $x^{\alpha}$, $x^{\beta}$, $x^{\gamma}$, etc entstehen und, indem man auch größere Zahlen als $N$ bei den Resten zulässt, werden aus den Potenzen $x^{2\alpha}$, $x^{3\alpha}$, $x^{4\alpha}$, etc die Reste $a^2$, $a^3$, $a^4$, etc entspringen, die also auch in der Reihe der Reste $1$, $a$, $b$, $c$, etc enthalten sein werden. Dann aber werden die Potenzen $x^{\alpha + \beta}$, $x^{\alpha + \gamma}$, $x^{\alpha + \beta + \gamma}$, etc die Reste $ab$, $ac$, $abc$, etc zurücklassen, die also auch in der Reihe der Reste der gefunden werden müssen werden. Die Produkte werden wie auch immer aus den Resten $1$, $a$, $b$, $c$, etc durch Multiplikation gebildet alle in derselben Reihe der Reste auftauchen, wenn natürlich die einzelnen durch Wegschaffen des Teilers $N$, sooft das gemacht werden kann, auf die kleinste Form gebracht werden.

\hfill\textsc{Q.E.D.}
\subsection*{Korollar 1}
\paragraph*{§43}
Diese Formen der Reste würden sich umso deutlicher zeigen, wenn anstelle derer jene Potenzen von $x$, woher sie entstanden sind, eingesetzt werden; dann tauchen nämlich natürlich nicht nur alle Potenzen dieser Potenzen, sondern auch irgendwelche Produkte derer in den Resten auf.
\subsection*{Korollar 2}
\paragraph*{§44}
Und dennoch wird daher die Anzahl der Reste nicht unbestimmt; so wie wir nämlich schon gesehen haben, dass aus unzähligen Potenzen die gleichen Reste hervorgehen, so werden, wenn all diese Reste, die aus gegenseitiger Multiplikation entstehen, auf die kleinste Form gebracht werden, auf eine gemäßigte Menge zurückgeführt.
\subsection*{Korollar 3}
\paragraph*{§45}
Wenn so die kleinste Potenz, die durch $N$ geteilt wiederum die Einheit zurücklässt, $x^{\nu}$ war, sodass die Anzahl der Reste $1$, $a$, $b$, $c$, etc gleich $\nu$ ist, dann werden in derselben Zahl alle Produkte, die aus der Multiplikation der Zahlen $a$, $b$, $c$, etc entstehen, enthalten sein, wenn natürlich von diesen der Teiler $N$ sooft, wie es gemacht werden kann, weggeschafft wird.
\subsection*{Bemerkung}
\paragraph*{§46}
Es wird ein einziges Beispiel genügen uns alle Zweifel, die vielleicht über diese sich zeigende Menge an Resten entstehen kann, aufzulösen. Es sei also $x=2$ und für den Teiler nehme man $N=15$, welcher natürlich zu $2$ prim ist; nun werden die einzelnen Potenzen von $2$ durch $15$ geteilt die folgende Reste zurücklassen:

\begin{tabular}{lllllllllllll}
Potenzen & $1$, & $2$, & $2^2$, & $2^3$, & $2^4$, & $2^5$, & $2^6$, & $2^7$, & $2^8$, & $2^9$, & $2^{10}$, & etc \\
Reste & $1$, & $2$, & $4$, & $8$, & $1$, & $1$, & $2$, & $4$, & $8$, & $1$, & $2$, & etc \\ 
\end{tabular}

Die Potenz also, die als erste wieder die Einheit ergibt, ist $2^4$, von welcher die Reste immer in derselben Reihenfolge $1$, $2$, $4$, $8$ wiederholt werden, sodass nur vier verschiedene Reste auftauchen. Hier ist nun klar, wie auch immer die Reste miteinander multipliziert werden, dass aber niemals Zahlen erzeugt werden, die nicht im selben Quadrupel eingeschlossen sind, nachdem sie natürlich durch Wegschaffen des Teilers $15$ auf die kleinste Zahl zurückgeführt worden sind. In diesem Beispiel tauchen auch unter den Resten nicht alle zu $15$ primen Teile auf, sondern es werden diese Teile $7$, $11$, $13$, $14$ ausgeschlossen, die in gleicher Weise zu $15$ prim sind; daher wird die oben [§34] gemachte Verteilung unter den zum Teiler primen Teile, die in den Resten auftauchten und die nicht auftauchten, illustriert, auf die man im folgenden besonders achte.
\section*{Theorem 9}
\paragraph*{§47}
Bei den Resten, die aus der Teilung der Potenzen einer Zahl durch einen zu ihr primen Teiler zurückbleiben, tauchen entweder alle zum Teiler primen Teile auf oder die Anzahl der nicht auftauchenden Teile wird gleich sein oder wird in einem vielfachen Verhältnis zur Zahl der Teile stehen, die die Reste festsetzen.
\subsection*{Beweis}
Es sei die Reihe der Potenzen $1$, $x$, $x^2$, $x^3$, $x^4$, $x^5$, etc und der zu $x$ prime Teiler sei $N$, dessen Anzahl der zu selbiger primen Teile gleich $n$ sei. Es sei weiter $x^{\nu}$ die kleinste Potenz, die durch $N$ geteilt wieder die Einheit zurücklässt, sodass die Anzahl aller verschiedenen Reste gleich $\nu$ ist; weil diese alle zu $N$ prime Zahlen sind, wird deren Anzahl entweder gleich $n$ oder kleiner sein, und im ersten Fall werden unter den Resten jedenfalls alle zu $N$ primen Teile auftauchen. Wir wollen also den Fall betrachten, in dem $\nu < n$ ist, und es seien	$1$, $a$, $b$, $c$, $d$, etc alle Reste, die aus der Teilung der Potenzen $1$, $x$, $x^2$, $x^3$, $x^4$, \dots $x^{n-1}$ durch den Teiler $N$ zurückgelassen werden; weil deren Anzahl gleich $\nu$ ist, werden dort nicht alle zu $N$ primen Teile auftauchen. Es sei also $\alpha$ der Teil dieser Art, der in den Resten nicht auftaucht und es kann bewiesen werden, dass auch keine dieser Zahlen $\alpha a$, $\alpha b$, $\alpha c$, $\alpha d$, etc in den Resten auftaucht. Denn wenn $\alpha a$ ein Rest wäre, der der Potenz $x^{\lambda}$ entspricht, würde, weil $a$ auch ein Rest ist, der aus einer Potenz, z.\,B. $x^{\xi}$, entsteht, $x^{\lambda} = an + \alpha a$ und $x^{\xi} = BN + a$ und daher $x^{\lambda} - \alpha x^{\xi} =  (A-\alpha B)N$ durch $N$ teilbar sein. Weil aber $x^{\xi}$ eine zu $N$ prime Zahl ist und $x^{\lambda} - \alpha x^{\xi} = x^{\xi}(x^{\lambda - \xi} - \alpha)$ ist, würde die Zahl $x^{\lambda - \xi} - \alpha$ durch $N$ teilbar und so würde die Potenz $x^{\lambda - \xi}$ durch $N$ geteilt entgegen der Annahme den Rest $\alpha$ übrig lassen. Weil also $\alpha$, $\alpha a$, $\alpha b$, $\alpha c$, etc, deren Anzahl gleich $\nu$ ist, zu $N$ prime Zahlen sind und durch Teilung durch $N$ auf zu $N$ prime Teile zurückgeführt werden können, wird man sofort auch einen einzigen Teil, der zu $N$ prim ist, in den Resten nicht finden, zugleich werden auch $\nu$ Teile solcher Art angegeben werden können, die in den Resten nicht auftauchen. Die Anzahl der nicht auftauchenden Teile ist also, wenn sie nicht $0$ ist, mindestens gleich $\nu$ und wenn außerdem ein zu $N$ primer Teiler $\beta$ in diesen Nicht-Resten nicht enthalten war, wird man erneut $\nu$ neue Teile haben, die in den Resten nicht auftauchen, und so weiter. Wenn daher nicht alle zum Teiler $N$ primen Teile in den Resten auftauchen, ist die Anzahl der Teile, die nicht auftauchen, notwendigerweise entweder gleich $\nu$ oder gleich $2\nu$ oder gleich $3\nu$ oder irgendein anderes Vielfaches von $\nu$, das heißt der Anzahl der verschieden Reste.

\hfill\textsc{Q.E.D.}
\subsection*{Korollar 1}
\paragraph*{§48}
Nachdem also der Unterschied zwischen den zum Teiler $N$ primen Teilen, die Reste sind, und denen, die keine Reste sind, festgesetzt worden ist, ist aus dem Beweis klar, dass das Produkt aus dem Rest und Nicht-Rest immer in der Klasse der Nicht-Reste enthalten ist. Wenn so $a$ ein Rest ist, $\alpha$ ein Nicht-Rest, wird also gewiss das Produkt $\alpha a$ kein Rest sein.
\subsection*{Korollar 2}
\paragraph*{§49}
Hingegen haben wir schon oben gesehen, dass das Produkt aus $2$ oder mehreren Resten in der Klasse der Reste gefunden wird. Daher folgt, dass das Produkt aus einem Nicht Rest und irgendwelchen Resten in der Klasse der Nicht-Reste auftauchen muss.
\subsection*{Bemerkung}
\paragraph*{§50}
Die Art dieses Beweis ist also auf dieses Fundament gestützt, wenn also unter den Resten diese zum Teiler primen Teile $1$, $a$, $b$, $c$, $d$, etc auftauchen und $\alpha$ auch ein zum Teiler primer Teil war, der in diesen Resten nicht enthalten ist, dass dann alle Produkte $\alpha a$, $\alpha b$, $\alpha c$, $\alpha d$, etc nicht nur in den Resten auftauchen, was freilich perfekt bewiesen worden ist, sondern dass sie auch zum Teiler $N$ prime Teile sind und alle untereinander verschieden sind oder, wenn sie tatsächlich durch $N$ geteilt werden, die verschiedenen Reste zurückgelassen werden. Jenes ist freilich per se klar; weil nämlich $\alpha$ wie $a$, $b$, $c$, $d$, etc zu $N$ prime Zahlen sind, ist notwendig, dass auch deren Produkte zu $N$ prim sind. Dass aber die Produkte $\alpha a$, $\alpha b$, $\alpha c$, $\alpha d$, etc, die auf $N$ bezogen, alle zueinander verschieden sind, sieht man ein, weil, wenn z.\,B. $\alpha a$ und $\alpha b$ durch $N$ geteilt die gleichen Reste gäben, deren Differenz $\alpha b - \alpha a = \alpha (b-a)$ durch $N$ teilbar wäre und daher auch $b-a$; das widerspricht aber der Annahme, dass $a$ und $b$ zu $N$ prime verschiedene Teile sind.
\section*{Theorem 10}
\paragraph*{§51}
Der Exponent der kleinsten Potenz $x^{\nu}$, die durch eine zu $x$ prime Zahl $N$ geteilt die Einheit übrig lassen, ist entweder gleich der Anzahl der zu $N$ primen Teile oder der Hälfte dieser Anzahl oder einem anderen echten Teil.
\subsection*{Beweis}
Es sei $n$ die Anzahl der zu $N$ primen Teile; weil $\nu$ von diesen die Reste festsetzen, wird die Anzahl der Nicht-Reste $n-\nu$ sein. Wir haben aber gesehen, dass diese Zahl entweder gleich $0$ oder gleich $\nu$ oder $2\nu$ oder ein anderes Vielfaches vom Exponent $\nu$ ist. Es sei also $n-\nu = (m-1)\nu$, sodass $m$ entweder die Einheit oder eine andere ganze Zahl bezeichnet, und daher werden wir $n=m\nu$ und $\nu = \frac{n}{m}$ erhalten; daher ist klar, dass der Exponent der kleinsten Potenz von $x$, die durch $N$ geteilt die Einheit zurücklässt, entweder gleich $\nu$ oder $m=1$ ist, oder gleich $\frac{n}{2}$, wenn $m=2$ ist, oder dass sie im Allgemeinen ein echter Teil der Zahl $n$ ist, die die Menge der zum Teiler $N$ primen Teile ausdrückt.

\hfill\textsc{Q.E.D.}
\subsection*{Korollar 1}
Wenn $x^{\nu}$ die kleinste Potenz war, die durch eine zu $x$ prime Zahl $N$ geteilt die Einheit zurücklässt, sind die folgenden Potenzen, die die selben Reste zurücklassen, $x^{2\nu}$, $x^{3\nu}$, $x^{4\nu}$, $x^{5\nu}$, etc und es sind zusätzlich keine anderen gegeben, die durch $N$ geteilt die Einheit zurücklassen.
\subsection*{Korollar 2}
\paragraph*{§53}
Der Exponent dieser kleinsten Potenz ist also immer mit der Anzahl der zu $N$ primen Teile so verknüpft, dass sie entweder jener selbst oder irgendeinem echten Teil von ihr gleich ist.
\subsection*{Bemerkung}
\paragraph*{§54}
Um dieses Verhältnis besser zu erkennen, wird es förderlich sein, einige einfache Fälle betrachtet zu haben. Es sei also $x=2$ und für $N$ wollen wir nacheinander ungerade zu $x=2$ prime Zahlen nehmen und wollen die kleinste Potenz von $2$ beschaffen, die durch eine ungerade Zahl geteilt die Einheit übrig lässt.

\begin{tabular}{>{\centering}p{1.7cm} >{\centering}p{4cm} >{\centering\arraybackslash}p{5cm}}
\toprule
\small\centering Teiler & \small \centering Anzahl der zu ihm primen Teile $n$ & \small\centering Die kleinste Potenz $2^{\nu}$, die durch $N$ geteilt die Einheit zurücklässt 
\tabularnewline
\midrule
\small
$3$ & $2$ & $2^2$, also $\nu = n$ \\
$5$ & $4$ & $2^{4}$, also $\nu = n$ \\
$7$ & $6$ & $2^{3}$, also $\nu = \tfrac{1}{2}n$ \\
$9$ & $6$ & $2^{6}$, also $\nu = n$ \\
$11$ & $10$ & $2^{10}$, also $\nu = n$ \\
$13$ & $12$ & $2^{12}$, also $\nu = n$ \\
$15$ & $8$ & $2^{4}$, also $\nu = \tfrac{1}{2}n$ \\
$17$ & $16$ & $2^{8}$, also $\nu = \tfrac{1}{2}n$ \\
$19$ & $18$ & $2^{18}$, also $\nu = n$ \\
$21$ & $12$ & $2^{6}$, also $\nu = \tfrac{1}{2}n$ \\
$23$ & $22$ & $2^{11}$, also $\nu = \tfrac{1}{2}n$ \\
$25$ & $20$ & $2^{20}$, also $\nu = n$ \\
$27$ & $18$ & $2^{18}$, also $\nu = n$ \\
$29$ & $28$ & $2^{28}$, also $\nu = n$ \\
$31$ & $30$ & $2^{5}$, also $\nu = \tfrac{1}{6}n$ \\

\bottomrule
\end{tabular}

\section*{Theorem 11}
\paragraph*{§55}
Wenn $N$ eine zu $x$ prime Zahl war und $n$ die Anzahl der zu $N$ primen Teiler, dann wird die Potenz $x^n$ um die Einheit vermindert immer durch die Zahl $N$ teilbar sein.
\subsection*{Beweis}
Es sei $x^{\nu}$ die kleinste Potenz, die durch $N$ geteilt die Einheit zurücklässt, und es wird $\nu$ entweder der Zahl $n$ gleich sein oder einem echten Teil $\frac{n}{m}$ von ihr. Weil also $x^{\nu}-1$ durch $N$ teilbar ist, wird, weil die Form $x^{\nu m}-1$ den Faktor $x^{\nu}-1$ hat, auch diese Form oder $x^n - 1$ durch $N$ teilbar sein.

\hfill\textsc{Q.E.D.}
\subsection*{Korollar 1}
\paragraph*{§56}
Wenn also der Teiler $N$ eine Primzahl ist und $x$ durch $p$ nicht teilbar ist, dann wird $x^{p-1}-1$ immer durch die Primzahl $p$ teilbar sein, wie ich schon längst bewiesen habe.
\subsection*{Korollar 2}
\paragraph*{§57}
Wenn außerdem $p$, $q$, $r$, etc Primzahlen sind und $x$ keine derer teilt, folgt aus diesem Theorem,
\begin{center}
\begin{tabular}{lc}
dass diese Formen &\quad durch \dots teilbar sein werden 
\tabularnewline
\small
$x^{p-1}-1$ & $p$ \\
$x^{(p-1)(q-1)}-1$ & $pq$ \\
$x^{pp(p-1)}-1$ & $p^3$ \\
$x^{p(p-1)(q-1)}-1$ & $ppq$ \\
$x^{(p-1)(q-1)(r-1)}-1$ & $pqr$ \\
\end{tabular}
\end{center}
\subsection*{Korollar 3}
\paragraph*{§58}
Wenn $x$ und $y$ zum Teiler $N$ prim sind, deren Anzahl der zu ihr primen Teile gleich $n$ sei, wird, weil $x^{n-1}$ und $y^{n}-1$ durch $N$ teilbar ist, auch $x^n - y^n$ immer durch die Zahl $N$ teilbar sein, welches Theorem allgemeiner ist.
\subsection*{Korollar 4}
\paragraph*{§59}
Nachdem also irgendeine Zahl $N$ vorgelegt wurde, deren Anzahl zu selbiger primen Teile gleich $n$ sei, wird, welche zu $N$ prime Zahl auch immer für $x$ genommen wird, die Formel $x^n = 1$ immer durch die Zahl $N$ teilbar sein.
\subsection*{Korollar 5}
\paragraph*{§60}
Oftmals kann es aber auch passieren, dass eine einfachere Form dieser Art, wie z.\,B. $x^{\frac{1}{2}n - 1}$ oder $x^{\frac{1}{3}n} - 1$ oder $x^{\frac{1}{4}}-1$, durch die Zahl $N$ teilbar ist, welcher Umstand von der Gestalt der Zahl $x$ abhängt.
\subsection*{Bemerkung}
\paragraph*{§61}
Sieh also diesen neuen Beweis des Fermat'schen Satzes, dass, wenn $p$ eine Primzahl war, alle Zahlen, die in dieser Form $a^{p-1}-1$ enthalten sind, durch $p$ teilbar sind, solange $a$ nicht durch $p$ teilbar ist. Ich hatte aber schon längst zwei Beweise dieses Theorems gegeben, aber der, den ich hier beschafft habe, scheint jenen vorzustehen, weil er nicht nur auf $p$ als Primzahl beschränkt ist. Denn welche Zahl $N$ auch immer für den Teiler angenommen wird, wird, solange $a$ zu ihr prim ist, diese Zahl $a^n - 1$ immer durch $N$ teilbar sein, wenn natürlich $n$ die Anzahl der zu $N$ primen Teile bezeichnet, welche Proposition sich um Vieles weiter erstreckt als die Fermat'sche. Daraus zeigt sich umso mehr die Nützlichkeit der ersten Theoreme, durch die ich die Anzahl der zu irgendeiner Zahl primen Teile bestimmt habe, was ohne diese Anwendung als zu unfruchtbar hätte erscheinen können.
\end{document}